\documentclass[12pt]{article}
\usepackage[german]{babel}
\usepackage{amssymb}
\usepackage{amsmath}
\usepackage{graphicx}
\usepackage{tikz,pgf}
\usetikzlibrary{calc}
\usetikzlibrary{decorations.pathmorphing}
\def\noi{\noindent}

\def\IN{\mathbb N}
\def\IZ{\mathbb Z}
\def\IR{\mathbb R}
\def\IC{\mathbb C}
\def\IQ{\mathbb Q}

\def\IK{\mathbb K}
\def\IL{\mathbb L}
\def\ik{\Bbbk}

\def\vs{\vspace}
\def\ma{\mathcal}

\begin{document}
	
	\begin{center}
		{\bf \Large Semialgebraicity of the convergence domain of an algebraic power series}
	\end{center}

	\vspace{0.5cm} \centerline{Tobias Kaiser}
	
	\vspace{0.7cm}
	\begin{center}
		\begin{minipage}[t]{10cm}\scriptsize{{\bf Abstract.}
				Given a power series in finitely many variables that is algebraic over the corresponding polynomial ring over a subfield of the reals, we show that its convergence domain is semialgebraic over the real closure of the subfield. This gives in particular that the convergence radius of a univariate Puiseux series that is algebraic in the above sense belongs to the real closure or is infinity.}
		\end{minipage}
	\end{center}
	
	\normalsize
	
	\section*{Introduction}
	
	The convergence radius of a univariate power series can be computed by the formula of Hadamard or, under some assumptions, by the formula of Euler.
	The formulas express the convergence radius as the limit (superior) of certain sequences stemming from the sequence of coefficients of the power series. Hence they are essentially of analytic and not of algebraic nature. We are interested in algebraic power series. These are the (necessarily convergent) power series that are algebraic over the polynomial ring over the field of complex numbers (see for example Wilczynski [12]). In the real case they occur as germs of Nash functions (see Bochnak, Coste and Roy [4]). More generally, we consider algebraic Puiseux series which appear as local descriptions of univariate algebraic functions (see for example Bliss [3]).
	
	The convergence radius is by definition a nonnegative real number or infinity. Bierstone [2] has dealt with the size of the convergence radius in the Nash case. We ask what can be said about the convergence radius if the annihilating polynomial has coefficients in a given subfield of the complex numbers. This setting has for example recently occurred in [8] where the periods from transcendental number theory have been studied in a geometric way. Of particular interest is the situation when the polynomial is over the algebraic numbers. And here the surprising fact (from the analytic point of view) is that indeed the convergence radius is a real algebraic number (or infinity). The same then holds for algebraic Puiseux series.

	As stated above, this result does not follow by computing the limit in Hadamard's or Euler's formula. The coefficients of an algebraic power or Puiseux series are determined by recursion formulas or closed formulas depending on the coefficients of the annihilating polynomial (see Hickel and Matusinski [5] and the literature cited there). But this does not seem to simplify the exact computation of the convergence radius.
	
	\enlargethispage{1cm}
	\rule{14cm}{0.01cm}
	
	{\footnotesize{\itshape 2020 Mathematics Subject Classification:} 14H05, 14P10, 14P15, 30B10,  32A05}
	\newline
	{\footnotesize{\itshape Keywords and phrases:} algebraic power series, algebraic Puiseux series, real closure, convergence radius, convergence domain}
	
    Our approach is of necessarily algebraic nature. We consider a subfield of the reals and a power series in a finite number of variables that is algebraic over the corresponding polynomial ring over the subfield.
	
    We are able to prove that its convergence domain is semialgebraic over the real closure of the subfield. This implies in the univariate case that the convergence radius of an algebraic power series or Puiseux series belongs to this real closure (or is infinity). We show that this result is optimal. In particular, we obtain the above result for the (real) algebraic numbers. The proof uses semialgebraic geometry as can be found in [4]  and basic facts about holomorphic functions in finitely many variables as can be found in H\"ormander [7]. Note that the results also hold for multivariate Puiseux series (see the setting of Hickel and Matusinski [6]) but we refrain from formulating these to keep the presentation short.
	
	\section*{Notations}
	
	By $\IN=\{1,2,3,\ldots\}$ we denote the set of natural numbers and by $\IN_0$ the set of natural numbers with $0$.
	
	For $a\in \IC$ and $r\in \IR_{>0}$ we set $B(a,r):=\{z\in \IC: |z-a|<r\}$. For $a\in \IC^n$ and $r\in\IR_{>0}^n$ we set
	$\Delta(a,r)=\Delta_n(a,r):=\prod_{i=1}^n B(a_i,r_i)$.
	In the case $a=0$ we often write $B(r)$ respectively $\Delta(r)$.
	
	For sets $X,Y$ we write $\pi_X:X\times Y\to X, (x,y)\to x,$ for the projection on the first component.
	Given a subset $A$ of $X\times Y$ and $x\in X$ we let $A_x:=\{y\in Y\mid (x,y)\in A\}$.

	\section*{Preliminaries}
	
    A power series $f(X)$ in the variables $X=(X_1,\ldots,X_n)$ with complex coefficients is written as $\sum_{\alpha\in \IN_0^n}a_\alpha X^\alpha$. It is convergent if there is $r\in \IR_{>0}^n$ such that $f(X)$ converges in $\Delta(r)$. 
	We denote by $$\ma{R}(f):=\big\{r\in \IR_{>0}^n\;\big\vert\; f\mbox{ converges on }\Delta(r)\big\}$$
	the set of {\bf relevant radii} of $f$. 
	The set $\ma{C}(f):=\bigcup_{r\in \ma{R}(f)}\Delta(r)$ is the {\bf convergence domain} of $f$. It is a logarithmically convex Reinhardt domain (see [7, Chapter 2.4]). We obtain a holomorphic function $f:\ma{C}(f)\to \IC, z\mapsto \sum_{\alpha\in \IN_0^n}a_\alpha z^\alpha$.
    In the univariate case we write $f(X)=\sum_{k=0}^\infty a_k X^k$. Then
    $R(f):=\sup\ma{R}(f)\in \IR_{\geq 0}\cup \{+\infty\}$ is the {\bf convergence radius} (where by definition $\sup\emptyset=0$).
    By Hadamard's formula $R(f)=1/\limsup_{k\to \infty}\sqrt[k]{|a_k|}$
    and by Euler's formula $R(f)=\lim_{k\to \infty}|a_{k+1}/a_k|$ if the latter exists. 
	We refer to Ruiz [10] for the general theory of power series.
	A Puiseux series is of the form $g(X)=\sum_{k=l}^\infty a_k X^{k/p}$ where $l\in \IZ$ and $p\in \IN$. Its convergence radius is given by $R(g)=R(f)^p$ where $f(X)=\sum_{k=0}^\infty a_kX^k$. Note that it converges on $B(R(g))\setminus \{0\}$ if there are nonvanishing coefficients for negative exponents. 
	We refer to Walker [11] for Puiseux series.

	\vs{0.5cm}
	Let $\ik$ be a subfield of the reals. We equip it with the ordering coming from the reals. By $\IK$ we denote the real closure of $\ik$ in $\IR$ and by $\IL$ its algebraic closure. Note that $\IL=\IK(i)$.
	We call a subset of $\IR^n$ or a function $\IR^n\to \IR$ $\IK$-semialgebraic if it is semialgebraic over $\IK$. We call a subset of $\IC^n$ $\IL$-Zariski closed or open if it is Zariski closed or open over $\IL$. We identify $\IC$ and $\IR^2$ canonically.
	
	\vs{0.5cm}
	A power series $f(X)$ is {\bf $\ik$-algebraic} if it is algebraic over the polynomial ring $\ik[X]$; i.e. there is $P(X,T)\in \ik[X,T]\setminus\{0\}$ with $P(X,f(X))=0$.
	We collect some basic facts and give references or proofs for the reader's convenience.
	
	\vs{0.5cm}
	{\bf Fact 1}
	
	\vs{0.1cm}
	Let $f(X)$ be a $\ik$-algebraic power series. 
	Then there is an irreducible $P(X,T)\in \ik[X,T]$ such that $P(X,f(X))=0$.
	
	{\bf Proof:}
	
	\vs{0.1cm}
	Let $Q\in \ik[X,T]\setminus\{0\}$ be with $Q(X,f(X))=0$. Let $Q_1(X,T),Q_2(X,T)\in \ik[X,T]$ be with $Q(X,T)=Q_1(X,T)Q_2(X,T)$.
	From
	$$0=Q(X,f(X))=Q_1(X,f(X))Q_2(X,f(X))$$
	we obtain $Q_1(X,f(X))=0$ or $Q_2(X,f(X))=0$ since the ring of power series is an integral domain.
	Hence we can pass to an irreducible factor of $Q$.
	\hfill$\blacksquare$
	
	\vs{0.5cm}
	{\bf Fact 2}
	
	\vs{0.1cm}
	The following holds:
	\begin{itemize}
		\item[(1)] A $\ik$-algebraic power series is convergent.
		\item[(2)] A $\ik$-algebraic power series has coefficients from $\IL$.
		\item[(3)] A power series is $\ik$-algebraic if and only if it is $\IK$-algebraic if and only if it is $\IL$-algebraic.
	\end{itemize}
	{\bf Proof:}
	
	\vs{0.1cm}
	Parts (1) and (2) follow from the Artin-Mazur description of algebraic power series (see [1] and [4, Chapter 8.1]).
	Part (3) follows from the transitivity of algebraicity (by Kaplansky [9, Remark after Theorem 22 in Chapter I \S 3] algebraicity over integral domains can be reduced to the field case where transitivity holds). 
	\hfill$\blacksquare$
	
	\vs{0.5cm}
	A $\ik$-algebraic Puiseux series is defined in the similar way. Facts 1 \& 2 hold similarly.
     
    \newpage
    \section*{Results}

    {\bf 1. Remark}
    
    \vs{0.1cm}
    Let $f$ be a convergent power series. Then $\ma{R}(f)$ is $\IK$-semialgebraic if and only if $\ma{C}(f)$ is $\IK$-semialgebraic.

    \vs{0.1cm} 
    {\bf Proof:}
     
    \vs{0.1cm}
    If $\ma{R}(f)$ is $\IK$-semialgebraic we get that
    $$\ma{C}(f)=\big\{z\in \IC^n\;\big\vert\; \exists\, r\in \ma{R}(f)\; \forall \,i\; |z_i|<r_i\big\}$$ is
    $\IK$-semialgebraic.
    If $\ma{C}(f)$ is $\IK$-semialgebraic we get that
    $$\ma{R}(f)=\big\{r\in \IR_{>0}^n\;\big\vert\;  \Delta(r)\subset \ma{C}(f)\big\}$$
    is $\IK$-semialgebraic.
    \hfill$\blacksquare$
     
    \vs{0.5cm}
    {\bf 2. Theorem}
     
    \vs{0.1cm}
    {\it The convergence domain of a $\ik$-algebraic power series is $\IK$-semialgebraic.}
     
    \vs{0.1cm}
    {\bf Proof:}
     
    \vs{0.1cm}
    Let $f$ be a $\ik$-algebraic power series.
    By Fact 2(3) we can assume that $\ik=\IK$.
    By Fact 1 there is an irreducible $P(X,T)\in \IK[X,T]$ with $P(X,f(X))=0$. 
    
     Let $A$ be the set of all $z\in \IC^n$ such that $P(z,t)=0$ implies $\big(\partial P/\partial T\big)(z,t)\neq 0$ for all $t\in \IC$.
     Then $A$ is $\IK$-semialgebraic. Its interior $B$ is also $\IK$-semialgebraic.
     
     \vs{0.2cm}
     {\bf Claim 1:}
     The set $B$ is dense in $\IC^n$. 
     
     \vs{0.1cm}
     {\bf Proof of Claim 1:}
     Let $Q:=\partial P/\partial T\in \IK[X,T]$.
     Since $P$ is irreducible in $\IK(X)[T]$ we obtain that $P$ and $Q$ are relatively prime in $\IK(X)[T]$.
     Hence we find $G,H\in \IK(X)[T]$ such that
     $GP+HQ=1$.
     Let $b\in \IK[X]$ be the common denominator of all the coefficients of $G$ and $H$.
     Then for every $z\in \IC^n$ with $b(z)\neq 0$ we have that there is no common zero of $P(z,T)$ and $Q(z,T)$.
     Hence $B\supset \{z\in \IC^n\mid b(z)\neq 0\}$. The latter is a nonempty $\IL$-Zariski open subset of $\IC^n$ and therefore dense in $\IC^n$.
     \hfill$\blacksquare_{\mathrm{Claim\,1}}$
     
     \vs{0.2cm}
     Let $c\in \IK[X]$ be the leading coefficient of $P$ with respect to the variable $T$. 
     We set $C:=\{z\in B\mid c(z)\neq 0\}$. Then $C$ is an open $\IK$-semialgebraic subset of $\IC^n$ that is dense.
     Let $d$ be the degree of $P$ with respect to $T$. Then for every $z\in C$ we have that $P(z,T)\in \IC[T]$ has degree $d$ and $d$ distinct zeros in $\IL$.
     For $z\in C$ let $N(z)$ be the set of zeros of $P(z,T)$.
     We set $N:=\{(z,t)\mid z\in C, t\in N(z)\}$.
     By the semialgebraic version of the implicit function theorem (see [4, Chapter 2.9]) we find a $\IK$-semialgebraic family
     $\Gamma\subset \IC^{n+1}\times \IR_{>0}^n\times \IC^{n+1}$
     such that $\pi_{\IC^{n+1}}(\Gamma)=N$ and the following holds.
     For every $(z,t)\in N$ there is a unique $s_{(z,t)}\in \IR_{>0}^n$ such that $(z,t,s_{(z,t)})\in \pi_{\IC^{n+1}\times \IR_{>0}^n}(\Gamma)$ and $\Gamma_{(z,t,s_{(z,t)})}$ is the graph of a holomorphic function $g_{(z,t)}:\Delta(z,s_{(z,t)})\to \IC$ with the following properties: 
     \begin{itemize}
     	\item[(1)] $\Delta(z,s_{(z,t)})\subset C$,
     	\item[(2)] $g_{(z,t)}(z)=t$,
     	\item[(3)] $P(w,g_{(z,t)}(w))=0$ for every $w\in \Delta(z,s_{(z,t)})$,
     	\item[(4)] $g_{(z,t_1)}(\Delta(z,s_{(z,t_1)}))\cap g_{(z,t_2)}(\Delta(z,s_{(z,t_2)}))=\emptyset$ for $t_1\neq t_2$. 
     \end{itemize}
     Note that the functions $g_{(z,t)}$ can be indeed chosen to be holomorphic by the classical inverse function theorem.
      
     Let $\ma{X}$ be a semialgebraic cell decomposition of $N$.
     For $X\in \ma{X}$ we set $Y_X:=\pi_{\IC^n}(X)$. 
     We have a continuous $\IK$-semialgebraic function $h_X:Y_X\to \IC$ such that $X=\mathrm{graph}(h_X)$. 
     We let $\ma{U}$ to be the set of all $X\in \ma{X}$ such that $Y_X$ is open in $\IC^n$.
     We set $\ma{V}:=\{Y_X\mid X\in \ma{U}\}$ and for $Y\in \ma{V}$ we define 
     $\ma{U}_Y:=\{X\in \ma{U}\mid Y_X=Y\}$. 
     
     \vs{0.2cm}
     {\bf Claim 2:} For $X\in \ma{U}$ we have that $h_X:Y_X\to \IC$ is holomorphic.
     
     \vs{0.1cm}
     {\bf Proof of Claim 2:} Let $Y:=Y_X$. By construction we have that $P(z,h_X(z))=0$ for all $z\in Y$.
     Let $z_0\in Y$ and $t_0:=h_X(z_0)$.
     By continuity we find $\rho\in \IR_{>0}^n$ with $\rho\leq s_{(z_0,t)}$ for every $t\in N(z_0)$ such that $h_X(\Delta(z_0,\rho))\cap h_{X'}(\Delta(z_0,\rho))=\emptyset$ for all $X'\in \ma{U}_Y\setminus\{X\}$. 
     This gives that $h_X$ coincides with $g_{(z_0,t_0)}$ on $\Delta(z_0,\rho)$ and hence is holomorphic there. Since $z_0$ is arbitrary we are done.
     \hfill$\blacksquare_{\mathrm{Claim\,2}}$

     \vs{0.5cm}
     For $r\in \IR_{>0}^n$ and $Y\in\ma{V}$ let
     $\mathrm{CC}(r,Y)$ be the set of connected components of $\Delta(r)\cap Y$.
     We set $\mathrm{CC}(r):=\bigcup_{Y\in \ma{V}}\mathrm{CC}(r,Y)$.
     For $W\in \mathrm{CC}(r)$ let $Y(W)$ be the unique element $Y$ of $\ma{V}$ with $W\subset Y$.
     Note that 
     the family 
     $$\big(\mathrm{CC}(r,Y)\big)_{r\in \IR_{>0}^n,Y\in \ma{V}}$$ 
     is $\IK$-semialgebraic. In particular there is a natural number $M$ such that 
     $$\# \mathrm{CC}(r)\leq M$$
     for all $r\in \IR_{>0}^n$.
     We call a map $\sigma:\mathrm{CC}(r)\to \ma{U}$ faithful if $Y_{\sigma(W)}=Y(W)$ for every $W\in \mathrm{CC}(r)$.
     By $\Sigma(r)$ we denote the set of all faithful maps $\sigma:\mathrm{CC}(r)\to \ma{U}$. 
     For $r\in \IR_{>0}^n$ let $W(r):=\bigcup_{W\in \mathrm{CC}(r)}W$. Note that $W(r)=\Delta(r)\cap \bigcup_{Y\in \ma{V}}Y$ and that $W(r)$ is dense in $\Delta(r)$.
     For $z\in W(r)$ let $W_z$ be the unique element of $\mathrm{CC}(r)$ where $z$ belongs to. For $\sigma\in \Sigma(r)$ we define
     $$\alpha_{r,\sigma}:W(r)\to \IC, z\mapsto h_{\sigma(W_z)}(z).$$
     Choose $r^*\in \ma{R}(f)$ such that $f|_{\Delta(r^*)}$ is $\IK$-semialgebraic. Note that such an $r^*$ certainly exists by [4, Corollary 8.1.6]. 
     We call $r\in \IR_{>0}^n$ good if there is $\sigma\in \Sigma(r)$ such that $\alpha_{r,\sigma}$ has a holomorphic extension $\beta_{r,\sigma}:\Delta(r)\to \IC$ with
     $$\beta_{r,\sigma}|_{\Delta(r)\cap \Delta(r^*)}=f|_{\Delta(r)\cap \Delta(r^*)}.$$
     The set of all $r\in \IR_{>0}^n$ which are good is $\IK$-semialgebraic. This follows since $W(r)$ is dense in $\Delta(r)$ for every $r\in \IR_{>0}^n$ and both continuous extendability and complex differentiabilty are first order properties. We show that $\ma{R}(f)=\{r\in \IR_{>0}^n\mid r\mbox{ good}\}$ and are done by Remark 1.
     
     Let $r\in \ma{R}(f)$.
     For $W\in \mathrm{CC}(r)$ choose $z_W\in W$. Let $t_W:=f(z_W)$. Then $t\in N(z_W)$. We obtain by the uniqueness of the implicit function that 
     $$f|_{\Delta(r)\cap\Delta(z_W, s_{(z_W,t_W)})}=g_{(z_W,t_W)}|_{\Delta(r)\cap\Delta(z_W, s_{(z_W,t_W)})}.$$
   
     By the identity theorem and the proof of Claim 2 we obtain that 
     $f|_W=h_{\sigma(W)}|_W$ where $\sigma(W)$ is the unique $X\in \ma{U}$ with $(z_W,t_W)\in X$. Note that $\sigma(W)$ does not depend on the choice of $z_W\in W$.
     The map $\sigma:\mathrm{CC}(r)\to \ma{U},W\to\sigma(W),$ is faithful.
     By construction we have that $\alpha_{r,\sigma}=f|_{W(r)}$.
     Hence $\alpha_{r,\sigma}$ has a holomorphic extension $\beta_{r,\sigma}:\Delta(r)\to \IC$ with
     $$\beta_{r,\sigma}|_{\Delta(r)\cap \Delta(r^*)}=f|_{\Delta(r)\cap \Delta(r^*)},$$
     namely $f|_{\Delta(r)}$.
    
     Let $r$ be good. Then there is $\sigma\in \Sigma(r)$ and a holomorphic function $\beta_{r,\sigma}:\Delta(r)\to \IC$ such that
      $$\beta_{r,\sigma}|_{\Delta(r)\cap \Delta(r^*)}=f|_{\Delta(r)\cap \Delta(r^*)}.$$ 

     Let $\tilde{f}$ be the power series expansion of $\beta_{r,\sigma}$ on $\Delta(r)$ which exists by [7, Theorem 2.4.5]. By the identity theorem we obtain that $\tilde{f}=f$.
     Hence $f$ converges on $\Delta(r)$ and we obtain that $r\in \ma{R}(f)$.  
     \hfill$\blacksquare$
     
     \vs{0.5cm}
     By [4, Corollary 8.1.6] an algebraic power series gives a semialgebraic function on some open neighbourhood of the origin. We obtain that the convergence domain can be taken.
     
     \vs{0.5cm}
     {\bf 3. Corollary}
     
     \vs{0.1cm}
     {\it Let $f(X)$ be a $\ik$-algebraic power series. Then the function
     	$f:\ma{C}(f)\to \IC, z\mapsto f(z),$ is $\IK$-semialgebraic.}
     
     \vs{0.1cm}
     {\bf Proof:}
     
     \vs{0.1cm}
     This has been shown in the proof of Theorem 2.
     \hfill$\blacksquare$

     \vs{0.5cm}
     The following corollary is trivial in the case $\IK=\IR$ but of interest in the other cases.

     \vs{0.5cm}
     {\bf 4. Corollary}
     
     \vs{0.1cm}
     {\it The convergence radius of a univariate $\ik$-algebraic Puiseux series belongs to $\IK_{>0}\cup\{\infty\}$.}
     
     \vs{0.1cm}
     {\bf Proof:}
     
     \vs{0.1cm}
     Let $g(X)=\sum_{k=l}^\infty a_k X^{k/p}$ be a $\ik$-algebraic Puiseux series where $l\in \IZ$ and $p\in \IN$. We have that $h(X)=\sum_{k=0}^\infty a_k X^{k/p}$ is also $\ik$-algebraic since this clearly holds for the difference $g(X)-h(X)$. 
     Let $Q\in \ik[X,T]\setminus \{0\}$ be such that $Q(X,h(X))=0$.
     Set $f(X):=\sum_{k=0}^\infty a_k X^k$ and $P(X,T):=Q(X^p,T)$.
     Then $f(X)=h(X^p)$ and $P(X,f(X))=0$. Hence $f$ is a $\ik$-algebraic power series.
     By Theorem 2 and Remark 1 we have that $\ma{R}(f)$ is $\IK$-semialgebraic. Since $R(f)=\sup\ma{R}(f)$ we get that $R(f)$ belongs to $\IK_{>0}\cup\{\infty\}$. Since $R(g)=R(h)=R(f)^p$ (where $\infty^p:=\infty$) we are done.
     \hfill$\blacksquare$

     \vs{0.5cm}
     Note that this is the best we can hope for.
     We denote by $\IR^\mathrm{alg}$ the field of real algebraic numbers. 
     
     \vs{0.5cm}
     {\bf 5. Example}
     
     \vs{0.1cm}
     Let $a\in \IR_{>0}^\mathrm{alg}$. Let $P(X,T):=(a-X)T-1$. 
     Then 
     $f(X):=1/(a-X)=\sum_{p=0}^\infty X^p/a^{p+1}$ has convergence radius $a$ and fulfils $P(X,f(X))=0$. Hence $f$ is $\IR^\mathrm{alg}$-algebraic. By Fact 2(3) it is $\IQ$-algebraic. 
     
     \vs{0.5cm}
     We denote by 
     $\IC^\mathrm{alg}$ the field of algebraic numbers.
     
     \vs{0.5cm}
     {\bf 6. Corollary}
     
     \vs{0.1cm}
     {\it The convergence domain of a $\IC^\mathrm{alg}$-algebraic power series is $\IR^\mathrm{alg}$-semialgebraic.}
     
     \vs{0.1cm}
     {\bf Proof:}
     
     \vs{0.1cm}
     We apply Theorem 2 to the case $\ik=\IK=\IR^{\mathrm{alg}}$ and $\IL=\mathrm{\IC}^\mathrm{alg}$ and use Fact 2(3).
     \hfill$\blacksquare$
     
     \vs{0.5cm}
     {\bf 7. Corollary}
     
     \vs{0.1cm}
     {\it The convergence radius of a univariate $\IC^\mathrm{alg}$-algebraic Puiseux series belongs to $\IR^\mathrm{alg}_{>0}\cup\{\infty\}$.}
     
     \vs{0.1cm}
     {\bf Proof:}
     
     \vs{0.1cm}
     We apply Corollary 4 to the case $\ik=\IK=\IR^{\mathrm{alg}}$ and $\IL=\mathrm{\IC}^\mathrm{alg}$ and use Fact 2(3).
     \hfill$\blacksquare$

	\vs{2cm}
	{\bf References}
	
	{\footnotesize
		\begin{itemize}
			
		\item[(1)] M. Artin, B. Mazur: On periodic points. 
		{Ann. of Math.} {\bf 81} (1965), 82-99.
		
		\item[(2)] E. Bierstone: Control of radii of convergence and extension of subanalytic functions.
		{\it Proc. Amer. Math. Soc.} {\bf 132} (2004), no. 4, 997-1003.

		\item[(3)] G. A. Bliss: Algebraic functions.
		Dover Publications, 1966.
			
		\item[(4)] J. Bochnak, M. Coste, M.-F. Roy: Real algebraic geometry. Ergebnisse der Mathematik und ihrer Grenzgebiete {\bf 36}, Springer, 1998.
		
		\item[(5)] M. Hickel, M. Matusinski:
		On the algebraicity of Puiseux series.
		{\it Rev. Mat. Complut.} {\bf 30} (3) (2017) 589-620.
		
		\item[(6)] M. Hickel, M. Matusinski:
		About algebraic Puiseux series in several variables.
		{\it J. Algebra} {\bf 527} (2019), 55-108.
			
		 \item[(7)] L. H\"ormander: An introduction to complex analysis in several variables.
	        North-Holland Publishing Co., 1990.
		
			\item[(8)] T. Kaiser: Periods, Power Series, and Integrated Algebraic Numbers.
			{\it Mathematische Annalen} (2024), DOI 10.1007/s00208-024-02802-2, 32 p.
			
			\item[(9)] I. Kaplansky:
			Commutative Rings. Allen and Bacon,  1970.
			
			\item[(10)] J. Ruiz: The Basic Theory of Power Series.
			Advanced Lectures in Mathematics. Vieweg, 1993.
			
			\item[(11)] R. J. Walker: Algebraic curves.
			Springer, 1978.
			
			\item[(12)] E. J. Wilczynski: On the Form of the Power Series for an Algebraic Function.
			{\it Amer. Math. Monthly} {\bf 26} (1919), no.1, 9-12.
	\end{itemize}}

	\vspace{1.5cm}
	\noi
	Tobias Kaiser \\
	Faculty of Computer Science and Mathematics\\
	University of Passau\\
	94030 Passau\\
	Germany\\
	email: tobias.kaiser@uni-passau.de
\end{document}